\documentclass[a4paper,11pt]{article}

\usepackage{latexsym,amssymb,amsmath,enumerate,geometry}
\usepackage[all]{xy}
\newtheorem{example}{Example}[section] 
\newtheorem{Def}[example]{Definition}  
\newtheorem{rem}[example]{Remark}  
\newtheorem{thm}[example]{Theorem} 
\newtheorem{alg}[example]{Algorithm} 
\newenvironment{proof}
  {\noindent {\bf Proof} } {\mbox{} \hfill $\Box$ \mbox{} \\ } 

\newcommand{\idX} {\lambda_{X}}
\newcommand{\idbA}{\lambda_{\bar{A}}}
\newcommand{\idbX}{\lambda_{\bar{X}}}
\newcommand{\idbY}{\lambda_{\bar{Y}}}

\newcommand{\grp} {\mathrm{grp}}
\newcommand{\mon} {\mathrm{mon}}
\newcommand{\sep} {\mathrm{sep}}
\newcommand{\sast}{_{\ast}}
\newcommand{\wtX} {\widetilde{X}}
\newcommand{\bbZ}{\mathbb{Z}}
\newcommand{\cA} {\mathcal{A}}
\newcommand{\cC} {\mathcal{C}}
\newcommand{\cF} {\mathcal{F}}
\newcommand{\cI} {\mathcal{I}}

\newcommand{\cL} {\mathcal{L}}
\newcommand{\cP} {\mathcal{P}}
\newcommand{\bcP}{\bar{\cP}}
\newcommand{\cQ} {\mathcal{Q}}
\newcommand{\cR} {\mathcal{R}}
\newcommand{\cRL}{\cR_{\cL}}
\newcommand{\cT} {\mathcal{T}}
\newcommand{\cY} {\mathcal{Y}}
\newcommand{\GAP} {{\sf GAP}}
\newcommand{\GAPt}{{\sf GAP3}}
\newcommand{\GAPf}{{\sf GAP4}}
\newcommand{\toR} {\rightarrow_{\cR}}

\newcommand{\eqR} {=_{\cR}}
\newcommand{\eqbR}  {=_{\bar{\cR}}}

\newcommand{\startoR} {\stackrel{*}{\rightarrow}_{\cR}}
\newcommand{\startoL} {\stackrel{*}{\rightarrow}_{\cL}}
\newcommand{\startoLQ}{\stackrel{*}{\rightarrow}_{\cL(\cQ)}}
\newcommand{\startoRL}{\stackrel{*}{\rightarrow}_{\cRL}}
\newcommand{\starfromtoR} {{\stackrel{*}{\leftrightarrow}}_{\cR}}
\newcommand{\barw}{\bar{\omega}}

\begin{document}

\begin{titlepage}
  \title
  { \Large \bf 
    Logged Rewriting Procedures with          \\
    Application to Identities Among Relations \\
    \mbox{}                                   \\
    \large
    \mbox{}                                   \\
    U.W.B. Mathematics Preprint 99.07         \\
    \mbox{}
  }
  \author
  { Anne Heyworth\thanks{Supported by a University of Wales,
                         Bangor, research studentship.}   \\
         and                                  \\
    Christopher D. Wensley                    \\ 
    {\small  School of Mathematics }          \\
    {\small  University of Wales, Bangor }    \\
    {\small  Gwynedd, LL57 1UT, U.K.}         \\ 
    \mbox{}                                   \\
    {\small  email: ~map130@bangor.ac.uk,
                    ~c.d.wensley@bangor.ac.uk} \\
    \mbox{}
  }

\end{titlepage}

\maketitle

\begin{center}

\mbox{}

rewrite system; Knuth-Bendix completion; group presentation; \\
crossed module; identity among relations.
\end{center}
\vspace{1cm}

%%%%%%%%%%%%%%%%%%%%%%%%%%%%%%%%%%%%%%%%%%%%%%%%%%%%%%%%%%%%%%%%%%%%%%%%%%%%%

\begin{abstract}  
We introduce logged rewrite systems and present a variation on 
the Knuth-Bendix algorithm for obtaining (where possible) 
complete logged rewrite systems.
This procedure is then applied to work of Brown and Razak Salleh, 
and an algorithm is developed which provides a set of generators 
for the module of identities among relations of a group presentation.
\end{abstract}
  
%%%%%%%%%%%%%%%%%%%%%%%%%%%%%%%%%%%%%%%%%%%%%%%%%%%%%%%%%%%%%%%%%%%%%%%%  
\pagebreak

\section{Introduction} \label{Intro}

In many mathematical calculations an object  $Q$  is reduced to an
equivalent object  $Q'$  by repeated application of operations
$T_1, T_2, \ldots$  taken from a suitable family of transformations.
A \emph{logged transformation} has the form  $Q \mapsto (L,Q')$
where the logged information  $L$  records the operations  $T_i$
in such a way that there is a function  $(L,Q') \mapsto Q$
recovering the original object.
For example, when  $Q$  is an element of a vector space with fixed basis,
and the  $T_i$  are invertible linear transformations,
we may take  $L$  to be the product of matrices  $M_i$  
representing the  $T_i$.
When  $Q, Q'$  are group presentations, and the  $T_i$  are
Tietze transformations, then  $L$  must record how new generators 
are introduced and old ones removed.

We are interested in the case of a group presentation
$\cP = \langle X ; R \rangle$  of a group  $G$
for which a complete rewrite system 
of the associated monoid presentation  $\bcP$  can be found.
By a refinement of the completion process we obtain a
complete logged rewrite system  $\cL(\cP)$  for  $G$.
Brown and Razak Salleh have given in Theorem 1.1 of \cite{BrSa}
a generating set  $\cI$  for  $\Pi_2(\bcP)$, the $\bbZ G$-module of 
\emph{identities among the relations} of  $\bcP$.
This construction requires a morphism  $k_1$  from the free groupoid 
on the Cayley graph of  $\cP$  to the free crossed module on  $\cP$.

The main results of this paper are Theorem \ref{sepdef2},
which shows how such a  $k_1$  may be obtained from the
logged information in  $\cL(\cP)$,
and Algorithms \ref{Reduce2}, \ref{KB2}
which implement these constructions.
Identities among the relations in  $\cP$  are found
by taking suitable words  $w$  such that there
is a logged reduction  $w \to (L_w,\lambda)$
where  $\lambda$  is the empty word and  $L_w$  
is one of the required identities.
As the process of Knuth-Bendix completion proceeds,
each critical pair determines either a new logged rule or an identity.

These two papers thus solve a sixty year old problem
of giving an algorithmic description of the second homotopy group
of a $2$-complex.

Once the set of identities  $\cI$  has been obtained,
the next stage is to determine a minimal generating set for  $\cI$
and to express each identity in terms of these generators.
Again a logged reduction can be performed,
and the resulting information used to obtain
\emph{identities among the identities}.

There are connections between this work and that of 
Squier, Lafont, Street, Groves, Cremanns and Otto on finite derivation types 
(see \cite{Lafont, Street, Groves, CrOt})
which will be discussed in a future paper \cite{BHW}.
For an approach to  $\Pi_2(\cP)$  using the method of pictures,
see for example Pride \cite{Pride}.

We would like to thank Ronnie Brown for many helpful discussions.

%%%%%%%%%%%%%%%%%%%%%%%%%%%%%%%%%%%%%%%%%%%%%%%%%%%%%%%%%%%%%%%%%%%%%%%

\section{Logged Rewrite Systems} \label{LRS}

For the basic notions of string rewriting systems we refer the
reader to \cite{BoOt} or \cite{CrOt}, and include here only
sufficient detail to fix our notation.

For  $A$  a set we denote by  $F(A)$  the free group on  $A$,
and by  $A^*$  the free monoid on  $A$,
both having the empty word  $\lambda_A$  as identity element.
We denote the length of a word  $w$  in  $F(A)$  or  $\bar{A}^*$  by  $|w|$.
A string-rewriting system  $\cR \subset A^* \times A^*$   
generates a reduction relation
$$
\toR \;=\; \{ (ulv, urv) : (l,r) \in \cR, ~u,v \in A^* \}
$$
on  $A^*$.
If  $(w_1,w_2) \in\, \toR$  then we write 
$w_1 \toR w_2$, and say that  $w_1$  reduces to $w_2$. 
The reflexive, transitive closure of  $\toR$  is written  $\startoR$, 
and the reflexive, symmetric, transitive closure
$\starfromtoR$  coincides with the congruence  $=_\cR$.
The monoid with presentation  $\mon \langle A, R \rangle$
is the factor monoid  $A^*/\!\!\eqR$.
We say that  $\cR$  is a \emph{rewrite system for a group}  $G$  on  $A$
if  $G \cong A^*/{\starfromtoR}$.

A \emph{well-ordering}  $>$  on  $A^*$  is an ordering 
in which there are no infinite sequences  $w_1 > w_2 > \cdots$  for 
$w_1, w_2, \ldots \in A^*$.
An ordering  $>$  on  $A^*$ is \emph{admissible} 
if  $y > z$  implies  $uyv > uzv$  for all  $u,v \in A^*$.
A rewrite system  $\cR$  is \emph{compatible} with an admissible well-ordering
if for all  $(l,r) \in \cR$  it is the case that  $l>r$.
An example of an admissible well-ordering is a 
\emph{length-lexicographic} ordering:
on choosing a total order on  $A$  we say that  $u>v$ 
if  $|u| > |v|$,  or if  $|u| = |v|$  and  $u$  is greater that  $v$  
in the lexicographic ordering induced by the order on $A$.

Two rules  $(l,r),\,(l',r')$  of a rewrite system  $\cR$  \emph{overlap}
if either (type 1)  $l$  is a subword of  $l'$, 
in which case we can write  $l' = ulv$;
or (type 2) the right-hand end of  $l'$  is equal to 
the left-hand end of  $l$, and we can write  $ul = l'v'$.  
In either case we can write  $ulv = l'v'$.
Applying the two rules, we obtain the \emph{critical pair} 
$(urv,r'v')$  resulting from the overlap. 
This critical pair can be \emph{resolved} if there exists  $w$  such that 
$urv \startoR w$  and  $r'v' \startoR w$.

If a rewrite system  $\cR$  is compatible with 
an admissible well-ordering  $>$, then  $\toR$  is \emph{Noetherian}. 
A rewrite system generating a Noetherian reduction relation is 
\emph{complete} if and only if all of its critical pairs can be resolved.
The \emph{Knuth-Bendix completion procedure} attempts to transform 
a given finite rewrite system into an equivalent complete one. 
It uses an admissible well-ordering on  $\bar{A}^*$ 
to orientate the rules of $\cR$, 
and then finds all the overlaps and computes the critical pairs. 
It attempts to resolve the critical pairs by reducing both sides 
with respect to  $\toR$. 
If the pair fails to resolve, 
then the reduced critical pair is added to  $\cR$
and the search for overlaps begins again. 
If the procedure terminates the rewrite system is complete.

For  $A$  a set, let
$\bar{A} ~=~ \{a^+  : a \in A\} \,\sqcup\, \{a^- : a \in A\}$
and define  ${\idbA}^- = \idbA$  and 
$$
(a^+)^+ = a^+, \;
(a^+)^- = a^-, \;
(a^-)^+ = a^-, \;
(a^-)^- = a^+
\quad \mbox{for all} \quad a \in A ~,
$$
so that  $w^+ = w$  and  $w^-$  are defined for all  $w \in \bar{A}^*$.

If  $R_A^{\pm} = \{(a^+a^-, \lambda_A) : \forall ~a \in \bar{A} \}$
there is a monoid isomorphism  
$$
\mu_A ~:~ F(A) \to \mon \langle \bar{A}, R_A^{\pm} \rangle,
          \; a \mapsto a^+, \; a^{-1} \mapsto a^- \quad
               \mbox{for all} \quad a \in A ~.
$$

Now let 
$\cP = \grp\langle X,\, \omega : R \to F(X) \rangle$ 
be a presentation of a group  $G$,
where the set  $R$  provides labels for the relators, 
and the function  $\omega$  identifies the relators 
as words in the free group.
We associate to  $\cP$  the monoid presentation
$\bcP \;=\; \mon \langle \bar{X},\,
    \barw ~:~ \bar{R} \to \bar{X}^* \rangle$,
where 
\begin{eqnarray*}
    \bar{R} & = & R \cup \bar{X} ~, \\
   \barw(x) & = & x^+ x^-, \quad x \in \bar{X}~, \\
\barw(\rho) & = & (\mu_X \circ \omega)\rho, \quad \rho \in R ~,
\end{eqnarray*}
so that  $\bar{X}^*/\!\!\eqbR \,\cong G$.
We may form an initial rewrite system
$\cR_{init} := \{ (\barw\rho, \idbX) : \rho \in \bar{R} \}$,
choose a total order on  $\bar{X}$,
and seek a complete rewrite system by applying the completion procedure
using the length-lexicographic order.

We now seek to extend  $\cR_{init}$  to a
\emph{logged rewrite system} for  $\bcP$.
This should be a set of rules which not only reduce any word in 
$\bar{X}^*$  to an irreducible word 
(unique if the rewrite system is complete)
but will express the actual reduction as a consequence of 
the original group relators.
In order to define this formally we recall the concept of crossed modules, 
as this is where \emph{consequences of the relators} 
are defined \cite{BrHu}.

Recall that a \emph{crossed $F$-module} is a pair
$\cC = (C,\delta)$  where  
$\delta : C \to F$  is a group homomorphism and  
$F$  acts (on the right) on  $C$  so that the following hold:
\begin{center}
\begin{tabular}{lll}
 CM1: &  $\delta(c^u)=u^{-1}(\delta c)u$
        & (pre-crossed $F$-module axiom), \\
 CM2: &  $c^{\delta c_1} = c_1^{-1}c c_1$
        &  (Peiffer relation). \\
\end{tabular}
\end{center}
If  $(C,\delta)$  and  $(D,\gamma)$  are crossed $F$-modules then a 
\emph{morphism of crossed $F$-modules} 
is a group homomorphism  $\theta : C \to D$  such that
$\theta(c^u) = (\theta c)^u$  and  
$\delta c = (\gamma \circ \theta)c$  for all  $c \in C, \,u \in F$.

Given a set  $R$,  a group  $F$  and a function $\omega : R \to F$, 
the \emph{free} crossed  $F$-module on $\omega$
is a crossed $F$-module  $(C,\delta)$ 
together with a function  $\omega_C : R \to C$
satisfying  $\delta \circ \omega_C = \omega$  if,
given any other crossed $F$-module  $(D,\gamma)$ 
with a map  $\omega_D : R \to D$
such that  $\gamma \circ \omega_D = \omega$,  
then there exists a unique morphism of crossed $F$-modules 
$\theta : C \to D$  which satisfies 
$\theta \circ \omega_C = \omega_D$.
Since  $\omega$  induces a homomorphism 
$\omega : F(R) \to F$,
this free crossed module  $(C,\delta)$ 
may be viewed as the induced crossed module
$\omega \sast (\cF)$  (see \cite{BrWe9}) where  $\cF$
is the crossed  $F(R)$-module  $(F(R),1)$
with conjugation action.
The universal property of the induced crossed module
states that the following diagram commutes:
$$ 
\xymatrix{
   &  &  D \ar[ddl]^{\gamma} \\
  F(R) \ar[urr]^{\omega_D} \ar[r]_(0.6){\omega_C} \ar[d]^{1}
   & C \ar[ur]|(0.4){\theta} \ar[d]_{\delta}  \\
  F(R) \ar[r]_(0.6){\omega}
   & F \\ 
}
$$  

The consequences of the relators  $R$  in the presentation  $\cP$
are elements of the free crossed  $F(X)$-module $(C(R),\delta_2)$ 
on the function  $\omega:R \to F(X)$.
As we shall be using rewriting techniques, we 
specify this crossed $F(X)$-module in terms of monoids and congruences.

If  $Y = R \times F(X)$  we write the elements of  $\bar{Y}$
as  $(\rho^{\epsilon})^u$, rather than  $(\rho,u)^{\epsilon}$,
where  $\epsilon \in \{+,-\}$.
Elements of the free monoid  $\bar{Y}^*$ 
are known as \emph{Y-sequences} and have the form  
$$
(\rho_1^{\epsilon_1})^{u_1} \cdots (\rho_n^{\epsilon_n})^{u_n} 
\quad \mbox{where} \quad 
\rho_i \in R, \, u_i \in F(X), \, \epsilon_i \in \{+,-\}.
$$ 

The action of  $F(X)$  on  $Y$, defined by 
$$
(\rho, u)^v = (\rho, uv)  \quad \mbox{for} \quad v \in F(X) ~,
$$
induces an action of  $F(X)$  on  $\bar{Y}^*$,
$$
((\rho^{\epsilon})^u)^v
 = (\rho^{\epsilon})^{uv} \quad \mbox{for} \quad v \in F(X) ~.
$$
We may also define an action of  $\bar{X}^*$  on  $\bar{Y}^*$  
generated by
$$
(\rho^\epsilon)^{x^+} = (\rho^\epsilon)^x, \quad
(\rho^\epsilon)^{x^-} = (\rho^\epsilon)^{x^{-1}},
$$
which agrees with the action of  $\mon \langle \bar{X}, R_X^{\pm} \rangle$
on  $\bar{Y}^*$  obtained from the  $F(X)$-action  via  $\mu_X^{-1}$.
Thus it is convenient to write terms in a Y-sequence as
$(\rho^\epsilon)^u$  rather than  $(\rho^\epsilon)^{\mu_X^{-1}(u)}$.

The monoid morphism  $\delta: \bar{Y}^* \to F(X)$ 
is defined to be that induced by  
$$
\delta( (\rho^+)^u ) ~ = ~ u^{-1} (\omega \rho) u ~,  \quad
\delta( (\rho^-)^u ) ~ = ~ u^{-1} (\omega \rho)^{-1} u  ~,
$$
and we let 
$\bar{\delta} = \mu_X \circ \delta ~:~ \bar{Y}^* \to \bar{X}^*$.
We may then define the Peiffer congruence  $=_P$  on  $\bar{Y}^*$ 
to be that generated by the set
$$  
P = \{(y^-zy^+, z^{\delta y}) : y,z \in \bar{Y} \}
    \; \cup \; \{(y^+ y^-,\idbY) : y \in \bar{Y} \} ~.
$$  
It can be verified that the free crossed $F(X)$-module on  
$\omega:R \to F(X)$  is  $(C(R), \delta_2)$  where
\begin{itemize}
\item  $C(R) = \bar{Y}^*/\!\!=_P$,
\item  $[c]$  denotes the class of  $c$  in  $C(R)$,
\item  the action of  $F(X)$  on  $C(R)$  is given by  $[c]^u = [c^u]$, 
\item  $\delta_2 [c] = \delta c$.
\end{itemize}
More detailed expositions of this construction may be found in
\cite{BrHu,Anne}.
The key idea is that `consequences of relators' are products of 
conjugates of relators which are represented by Y-sequences. The notion
of equality of consequences of relators is expressed by the Peiffer
relations on the Y-sequences.

\begin{Def}  
A \emph{logged rewrite system} for the group presentation  $\cP$  of  $G$ 
is a set of triples 
$$
\cL = \{ (l_1, c_1, r_1), \ldots, (l_n,c_n,r_n) \} ~,
$$  
where  $c_1, \ldots, c_n \in \bar{Y}^*$
and    $\cRL = \{ (l_1,r_1), \ldots, (l_n,r_n)\}$ 
is a rewrite system for  $G$  on  $\bar{X}$ 
such that  $l_i = (\bar{\delta} c_i) r_i$  for  $i=1, \ldots, n$.
The system  $\cL$  is \emph{complete} if  $\cRL$  is complete.  
The \emph{initial logged rewrite system} of  $\cP$  is
$$
\cL_{init} \; = \;
  \{~(\barw\rho,\, (\rho)^{\idX},\, \idbX) ~:~ \rho \in R\} 
  ~\cup~
  \{~(x^+x^-, \idbY, \idbX) ~:~ x \in \bar{X} \} ~.
$$
\end{Def}  

\begin{alg}[Logged rewriting]  \label{Reduce2}
Given a logged rewrite system  $\cL$  and a word  $w \in \bar{X}^*$,
an irreducible word  $I(w)$  is determined, together with a Y-sequence  
$L(w)$,  such that  $w \startoRL I(w)$ 
and  $w = \delta_2(L(w)) I(w)$.
\begin{enumerate}[R1:]  
\item
(Initialise) Set $i = 0, \; z_i = w, \; L_i(w) = \idbY$. 
Clearly $w = \delta_2(L_i(w)) z_i$.
\item
(Search) 
Find  $(l,c,r) \in \cL$  and  $u,v \in \bar{X}^*$  such that  $z_i = ulv$ 
(i.e.  $z_i$  is reducible with respect to  $\cRL$).
\item
(Rewrite and Record)
When  $(l,c,r),u,v$  are found, increment  $i$,
set  $z_i = urv$  and set  
$L_i(w)$ $\,=\,$ $L_{i-1}(w) c^{u^{-1}}$.
It is easily verified at each stage that $w = \delta_2(L_i(w)) z_i$.
\item
(Loop) Repeat Search(R2:) and Rewrite(R3:) until Search fails.
\item
(Terminate) Output  $I(w) = z_i$  and  $L(w) = L_i(w)$.
\end{enumerate}
\end{alg}  
  
We represent this logged reduction as:
``$w \startoL I(w) \mbox{ by }  L(w)$''.
Note that if  $\cL$  is complete then there is a unique irreducible
word  $I(w)$  for each equivalence class  $[w]$  
under the congruence generated by  $\cL$.
Thus we obtain a section  $\sigma : G \to F(X)$  of  $\phi$,
mapping each group element to its unique representative in  $F(X)$,
and a \emph{normal form} function
$N = \sigma \circ \phi : F(X) \to F(X)$
with  $(N \circ \sigma)g = \sigma g$  for all  $g \in G$.

We now extend the process of completing a rewrite system,
which adds unresolved critical pairs resulting from overlaps,
by recording the logged part of each rule used.
This logged Knuth-Bendix procedure results in a rewrite system which stores 
information on how the rules were constructed.
The information is not unique because a rule may be derived 
in many different ways.

\begin{alg}[Logged Knuth-Bendix Completion] \label{KB2}
Given a group  $G$  with presentation 
$\cP=\grp\langle X, \omega:R \to F(X) \rangle$,
and an ordering on  $\bar{X}^*$  such that the standard Knuth-Bendix 
algorithm completes, a complete logged rewrite system  $\cL(\cP)$ 
is produced which respects the chosen ordering,
together with a set  $\cI(\cP)$  of identities among the relations.
\begin{enumerate}[K1:]  
\item
(Initialise)
Set  $j = 0$  and set  $\cL_j$
to be the logged rewrite system  $\cL_{init}(\cP)$.\\
Denote by  $I_j(w), \, L_j(w)$  the output from Algorithm \ref{Reduce2} 
using  $\cL_j$  when  $w$  is input.
\item
(Search for critical pairs)
Set  $crit = \emptyset$  and  $ident = \emptyset$.\\
Find  $(l,c,r),\,(l',c',r') \in \cL_j$
and  $u, \, v, \, v' \in \bar{X}^*$
such that  $u l v = l' v'$  where
(type 1)  $v' = \idbX$,  or 
(type 2)  $v = \idbX$. \\
Put  $
z = I_{j-1} (u r v), \;
d = L_{j-1} (u r v), \;
z' = I_{j-1} (r' v')$  and  $
d' = L_{j-1} (r' v')$.\\
If  $z \ne z'$  then add the logged critical pair
$(z', ({d'}^-) ({c'}^-) (c^+)^{u^{-1}} (d^+), z)$  to  $crit$,\\
otherwise add the Y-sequence  
$({d'}^-) ({c'}^-) (c^+)^{u^{-1}} (d^+)$  to  $ident$.
\item
(Add logged rules) 
Set  $\cL_j = \cL_{j-1}$.
For each logged critical pair  $(z',(c),z)$  in  $crit$ 
\begin{itemize}
\item add the logged rule
$(z', (c^+), z)$  to  $\cL_j$  if  $z < z'$,  or
\item add the logged rule
$(z, (c^-), z')$  to  $\cL_j$  if  $z > z'$.
\end{itemize}
\item
(Remove redundant rules)
For each new rule  $\psi \,\in\, \cL_j \setminus \cL_{j-1}$,  
test whether  $\psi$  reduces to  $(\idbX,\idbY,\idbX)$
using the remaining rules.
If so, remove  $\psi$  from  $\cL_j$.
\item
(Loop)
Repeat: increment  $j$, 
Search(K2:) for unresolvable critical pairs,
Add(K3:) logged rules, and
Remove(K4:) redundant rules
until all critical pairs are resolved.
\item
(Terminate)
Output the final logged rewrite system $\cL(\cP) = \cL_j$,
together with the list of identities  $\cI(\cP) = ident$.
\end{enumerate}
\end{alg}
 
It can be seen immediately from the description of the algorithm that,
if the middle term of each triple is omitted,
the standard completion of the rewrite system is obtained.
Also, for each  $j$,  the triples  $(l,c,r)$  which are added to  $\cL_j$
satisfy the requirement  $l = (\delta_2 c) r$.  
So, if the completion procedure terminates,  
$\cL(\cP)$  is a complete logged rewrite system.
The above procedures are illustrated in the following example.

Recall that the \emph{root} of a word  $w \in \bar{X}^*$
is the smallest initial subword  $v$  of  $w$  such that
$w = v^m$  for some  $m \in \mathbb{Z}^+$.
When  $w = \barw \rho$  there are overlaps in the words  
$v^{-1} v^m$  and   $v^m v^{-1}$, 
giving two logged critical pairs  
$(v^{m-1}, (\rho^+)^v, v^{-1})$  and  $(v^{m-1}, (\rho^+), v^{-1})$,
so we obtain the \emph{root module identity}
$\iota_{\rho} = (\rho^+)^v (\rho^-)$.
Note that, if  $w = v$, then  $[\iota_{\rho}] = [\lambda_{\bar{Y}}]$.

\begin{rem} \label{simprules}
The Y-sequence  $L(w)$  is constructed as an element of the 
free pre-crossed module on  $\omega$. 
Since we are working in the crossed module  $C(R)$
we may choose a simpler representative from  $[L(w)]$,  
using the Peiffer rules or the root module identities.
So if   $(l,c,r) \in \cL(\cP)$  and  $c$  contains  $c_1$  as a subsequence,
then applying one of the \emph{simplification rules}
\begin{enumerate}[S1:]
\item  $(\rho^{\epsilon})^v  \;\;\mbox{and}\;\;
        (\rho^{\epsilon})^{v^{-1}}  \;\to\;
        (\rho^{\epsilon})$,
       where  $v$  is the root of  $\omega \rho$, 
\item  $(\rho_1^{\epsilon_1})^{u_1} \, (\rho_2^{\epsilon_2})^{u_2}  \;\to\;
        (\rho_2^{\epsilon_2})^{u_2}
         (\rho_1^{\epsilon_1})^{u_1 u_2^{-1}
          (\omega \rho_2)^{\epsilon_2} u_2}   \;\;\mbox{or}\;\;
        (\rho_2^{\epsilon_2})^{u_2 u_1^{-1}
         (\omega \rho_1)^{-\epsilon_1} u_1} \,
          (\rho_1^{\epsilon_1})^{u_1}$, 
\end{enumerate}
to  $c_1$  gives an alternative logged rewrite system for  $\cP$.
\end{rem}

\begin{rem} \label{primidy}
The identity list  $\cI(\cP)$  is obtained as a by-product of the algorithm.
It may be simplified by deleting duplicates and conjugates; by applying the
two simplification rules; and by searching to see if one identity occurs
as a subsequence of another.

Any identity which satisfies the \emph{primary identity property}
(see \cite{BrHu}) is equivalent to the trivial identity, 
and may be omitted from the list.
An identity  
$(\rho_1^{\epsilon_1})^{u_1}(\rho_2^{\epsilon_2})^{u_2}
 \cdots (\rho_r^{\epsilon_r})^{u_r}$
satisfies this property if the set  $\{1,\ldots,r\}$
can be partitioned into pairs  $(i,j)$  such that for each pair
$\rho_i = \rho_j, \; \epsilon_i = - \epsilon_j$,  and
$u_i u_j^{-1}$  lies in the normal closure of the subgroup of
$F(X)$  generated by the relators  $\omega R$.
\end{rem}

\begin{example} \label{exquat}
\emph{The quaternion group  $Q_8$  is presented by 
$\cQ = \grp \langle X,\, \omega : R \to F(X) \rangle$  where
}  
$$
X = \{a,b\}, \;\; 
R = \{\rho_1, \rho_2, \rho_3, \rho_4\}, \;\;
\omega : \rho_1 \mapsto a^4, \;       \rho_2 \mapsto b^4, \;
         \rho_3 \mapsto abab^{-1}, \; \rho_4 \mapsto a^2b^2 ~.
$$
\emph{We begin with the logged rewrite system 
\begin{align*}
\cL_{init}(\cQ) \, = \,
  &  \{(a^{+4}, (\rho_1^+), \idbX),
     \ (b^{+4}, (\rho_2^+), \idbX), 
     \ (a^+b^+a^+b^-, (\rho_3^+), \idbX),\\
  &  \ (a^{+2}b^{+2}, (\rho_4^+), \idbX),
     \ (a^+a^-, \idbY, \idbX),
     \ (a^-a^+, \idbY, \idbX),\\
  &  \ (b^+b^-, \idbY, \idbX),
     \ (b^-b^+, \idbY, \idbX) \}.
\end{align*}
\noindent
We look for overlaps between the left hand sides of the rules. 
The first overlap is in the word  $a^-a^{+4}$,  which contains the
two left hand sides  $a^-a^+, \, a^{+4}$,  giving the rule
$(a^{+3}, (\rho_1^+)^a, a^-)$.  
As we have seem in Theorem \ref{simprules},
the logged part may be simplified to  $(\rho_1^+)$,  
since we could have used overlaps in  $a^{+4}a^-$.
The next three reductions, on words
}
$$
a^-a^+b^+a^+b^-, \quad 
a^-a^{+2}b^{+2}, \quad
\mbox{and} \quad
b^-b^{+4},
$$
\emph{give rules
$$
(b^+a^+b^-, (\rho_3^+)^a, a^-), \quad
(a^+b^{+2}, (\rho_4^+)^a, a^-), \quad 
\mbox{and} \quad
(b^{+3}, (\rho_2^+)^b, b^-)$$
respectively.
The next overlap is the first with two non-inverse rules: 
}
$$  
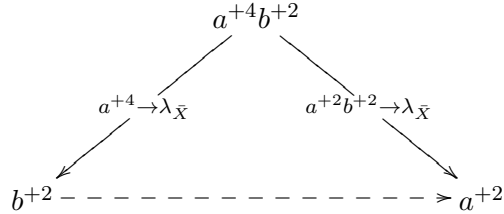
\xymatrix{
    & & a^{+4}b^{+2} \ar[ddll]|{a^{+4} \to \idbX}
                      \ar[ddrr]|{a^{+2}b^{+2} \to \idbX}  & &  \\
    & & & &  \\  
  b^{+2} \ar@{-->}[rrrr] & & & & a^{+2}  \\}  
$$  
\emph{Without the log part the critical pair is  $(b^{+2}, a^{+2})$. 
For the logged rewrite rule we need  $[c] \in C(R)$ 
so that  $b^{+2} = \delta_2(c) a^{+2}$  
where  $c$  is a product of conjugates of relators. 
The new logged rewrite rule as defined in step K3 is
$(b^{+2}, (\rho_1^-)(\rho_4^+)^{a^{-2}}, a^{+2})$ 
so  $c = (\rho_1^-) \,(\rho_4^+)^{a^{-2}}$.
This is verified by:
}
$$
\delta_2((\rho_1^-)(\rho_4^+)^{a^{-2}})\,a^{+2}
 \; = \;  (a^{-4})\, (a^{+2} (a^{+2}b^{+2}) a^{-2})\, a^{+2}
 \; = \;  b^{+2} ~.
$$  
\emph{Continuing with the logged completion for $\cQ$, a total of
$44$ rules are formed, of which  $22$  become redundant.
A complete logged rewrite system  $\cL(\cQ)$  is obtained
containing the list of rules  $(l_i,c_i,r_i)$ 
shown in Table 1.  %% \ref{lrwsq}.
}

\begin{table}[!htp] \label{lrwsq}
\begin{center}
\begin{tabular}{|r|l|l|l|}
\hline
$i$  &  $l_i$  &  $c_i$  &  $r_i$  \\
\hline
1 & $a^+a^-$  & $\idbY$ & $\idbX$ \\
2 & $a^-a^+$  & $\idbY$ & $\idbX$ \\
3 & $b^+b^-$  & $\idbY$ & $\idbX$ \\
4 & $b^-b^+$  & $\idbY$ & $\idbX$ \\
5 & $b^+a^+$
& $ (\rho_3^+)^{a^+} \,(\rho_1^-) \,(\rho_4^+)^{a^-}$
& $ a^+b^-$ \\
6 & $b^{+2} $
& $  (\rho_1^-) \,(\rho_4^+)^{a^{-2}} $
& $  a^{+2} $  \\ 
7 & $b^+a^-$ 
& $  (\rho_3^-)$
& $  a^+b^+$ \\
8 & $a^-b^+$
& $ (\rho_1^-) \,(\rho_4^+)^{a^-}$
& $ a^+b^-$  \\
9 & $a^{-2} $
& $  (\rho_1^-)$
& $ a^{+2}$  \\
10 & $a^-b^-$
& $  (\rho_1^-) \,(\rho_2^-)^{a^{-2}} \,(\rho_4^+)
     \,(\rho_3^+)^{a^+b^-} \,(\rho_3^-)$
& $ a^+b^+$  \\
11 & $b^-a^+ $
& $  (\rho_4^-) \,(\rho_3^+)^{a^-}$
& $ a^+b^+$  \\
12 & $b^-a^- $
& $  (\rho_3^-)^{a^+b^+}$
& $ a^+b^-$  \\
13 & $b^{-2}$
& $  (\rho_4^-)$
& $  a^{+2}$  \\ 
14 & $a^{+3} $
& $  (\rho_1^+)$
& $ a^-$  \\
15 & $a^{+2}b^+$
& $  (\rho_4^+)$
& $ b^-$  \\ 
16 & $a^{+2}b^- $
& $  (\rho_4^-)^{a^{-2}} \,(\rho_1^+)$
& $b^+$ \\
\hline
\end{tabular}
\end{center}
\caption{Complete logged rewrite system  $\cL(\cQ)$}
\end{table}

\noindent
\emph{So, for example, $a^+b^+b^+a^+$  reduces to  $\idbX$  as follows:
}
$$
a^+b^+b^+a^+
  \;\startoLQ\; a^{+4} \; (\mbox{by} \;  {c_{6}}^{a^-})
  \;\startoLQ\; a^+a^- \; (\mbox{by} \; {c_{14}}^{a^-})
  \;\startoLQ\; \idbX  \; (\mbox{by} \;  c_{1})
$$
\emph{and the logged information is given by:
}
$$
L(a^+b^+b^+a^+) 
  \;  =  \; (\rho_1^-)^{a^-}(\rho_4^+)^{a^{-3}}(\rho_1^+)
  \; \to \; (\rho_4^+)^{a^+}  \; \mbox{by (S1), (S2)} ~.
$$
\end{example}

%%%%%%%%%%%%%%%%%%%%%%%%%%%%%%%%%%%%%%%%%%%%%%%%%%%%%%%%%%%%%%%%%%%%%%%%%%%%%

\section{Computing a Set of Generators for $\Pi_2$}  

We now apply our procedures to the work of Brown and Razak Salleh 
in \cite{BrSa}.
Consider the short exact sequence
$$
\xymatrix{C(R) \ar[r]^{\delta_2} & F(X) \ar[r]^\phi & G
\ar[r] & 1}
$$ 
where  $\phi : F(X) \to G$  is the factor morphism.
By exactness, any element of  $F(X)$  that represents the identity in  $G$ 
is expressible as a (non-unique) product of conjugates of relators. 
Logged rewrite systems provide a method for doing this, 
and we have seen that a complete rewrite system for  $G$  
defines a section  $\sigma$  of  $\phi$
and a normal form function  $N$  on  $F(X)$.
The kernel $\ker \delta_2 = \Pi_2(\cP)$, 
the $\bbZ G$-module of identities among relations,
is of particular interest and logged rewrite systems help in the
calculation of a set of generators for it.
The reader is referred to \cite{BrSa} for details of contracting
homotopies and crossed resolutions. 
Here we merely use the formulae and results of that paper.
The main result we use is that a complete set of generators for the module
of identities among relations may be obtained from the set of relator cycles
of the Cayley graph, by defining particular functions and applying them
to the relator cycles.

The Cayley graph  $\wtX$  of the presentation
$\cP = \grp \langle X,\, \omega:R \to F(X) \rangle$ 
has the elements of  $G$  as its vertices.  
Edges of  $\wtX$  are written as pairs  $[g,x]$, 
where  $g$  is the group element identified with the source vertex, 
$x$  is a group generator identified with the edge label,
and  $g(\phi x)$  is the target vertex.
The crossed  $F(X)$-module  $(C(R), \delta_2)$  
is as defined in Section \ref{LRS}.
The free groupoid of all paths on  $\wtX$  is denoted
$F(\wtX)$ 
and paths are written  $[g,u]:g \to g (\phi u)$ 
for  $g \in G, \; u \in F(X)$.
We now quote Theorem 1.1 of  \cite{BrSa} which 
defines a generating set of identities among relations for $\cP$
as a  $\bbZ G$-module. 

\begin{thm} \label{sepdef1}
The module  $\Pi_2(\cP)$  is generated as a $\bbZ G$-module by elements
$$
\sep[g,\rho] = (\rho^-)^{(\sigma g)^{-1}} \, (k_1[g,\omega \rho])
$$
for all  $g \in G, \, \rho \in R$, where 
\begin{enumerate}[i)]
\item
$\sigma : G \to F(X)$  is a section of the quotient morphism $\phi:F(X)
\to G$,
\item
$k_1$  is a groupoid morphism from  $F(\wtX)$  to the one object
groupoid  $C(R)$,
\item
$\delta_2 k_1[g,x]
  =  (\sigma g)\,x\,(\sigma(g(\phi x)))^{-1}$ for all $x \in X, \, g \in G$.
\end{enumerate}
\end{thm}

The identities  $\sep[g,\rho]$  may be seen as \emph{separation elements} 
in the geometry of the Cayley graph with relators. 
We note that this complete set of generators is usually not minimal 
but we are not concerned with that here.
The main point of this part of our paper is to show that a logged
rewrite system for a presentation provides constructions of the functions
$\sigma$  and  $k_1$  satisfying the conditions given in the
theorem. 
Thus the results of \cite{BrSa} can be used, together with
logged rewriting procedures, to specify a generating set for $\Pi_2$ as
a $\bbZ G$-module. 
Implementation of the resulting procedure is discussed in Section \ref{imp}.
Minimising the generating sets thus obtained is a problem resolved in
the sequel to this paper \cite{HeRe}.

\begin{thm}[Separation Morphism] \label{sepdef2}
If the group presentation  $\cP$
has a complete logged rewrite system  $\cL(\cP)$
then the logged information determines a morphism  
$k_1 : F(\wtX) \to C(R)$ 
such that the elements  $\sep[g,\rho]$   
for all  $g \in G$  and  $\rho \in R$ 
form a complete set of generators for  $\Pi_2(\cP)$  
as a  $\bbZ G$-module.   
\end{thm}

\begin{proof}
For  $[g,x] \in \wtX$  choose an element in 
$[L((\sigma g) x (\sigma(g(\phi x)))^{-1})]$
as  $k_1[g,x]$,  choosing  $\idbY$  whenever possible.
These choices induce a groupoid morphism  $k_1 : F(\wtX) \to C(R)$,   
and  $\delta_2 k_1[g,x] = $  $(\sigma g) x (\sigma(g(\phi x)))^{-1} $.
Therefore, by the previous Theorem, the elements 
$\sep[g, \rho]$  for  $g \in G, \, \rho \in R$ 
are a complete set of generators for  $\Pi_2(\cP)$.
\end{proof}

%%%%%%%%%%%%%%%%%%%%%%%%%%%%%%%%%%%%%%%%%%%%%%%%%%%%%%%%%%%%%%%%%%%%%%%%%%%%

\section{Implementation} \label{imp}

A collection of functions is included in the first author's thesis \cite{Anne},
written using the computational group theory program {\GAPt}, 
to perform these calculations.
These functions have been rewritten for {\GAPf} \cite{GAP}
and submitted as a share package {\sf IDREL}.  
The structure of the program is outlined in the following algorithm.

\begin{alg}[Identities Among Relations]
Given a presentation  $\cP = \grp \langle X, \omega:R \to F(X) \rangle$ 
of a finite group $G$, a set of Y-sequences is determined
whose Peiffer equivalence classes generate  $\Pi_2(\cP)$  
as a  $\bbZ G$-module.
\begin{enumerate}[B1:]
\item
(Logged Rewrite System)
Apply Algorithm \ref{KB2} to obtain the completion  $\cL(\cP)$  from
the logged rewrite system  $\cL_{init}(\cP)$.
Let  $L$  and  $N$  be the log and normal form functions  
determined by Algorithm \ref{Reduce2} using  $\cL(\cP)$.
\item
(Cayley Graph)
The Cayley graph is represented by a list of edges, which are pairs  
$[g,x]$  where  $g$  is an irreducible word 
(with respect to  $\cL(\cP)$)
in  $F(X)$  and  $x \in X$.  
\item
(Definition of $k_1$)
The map  $k_1$  is defined on the edges by
$k_1[g,x]= \idbY$  if  $gx$  is irreducible with respect to  $\cL(\cP)$, and  
$k_1[g,x] \in [L((\sigma g) x (\sigma(g(\phi x)))^{-1})]$ otherwise.
\item
(Determination of Identities)
All pairs  $[g,\rho]$ 
where  $g$  is a vertex and  $\rho$  is a relator are considered.  
The boundary of the cycle is found by splitting up the  
relator  $\omega \rho$  to obtain a list of edges.  
The remaining edges of each  
cycle are identified with their images under  $k_1$.  
The identities are calculated by manipulating the information held so as  
to obtain a Y-sequence representing 
$\iota_{[g,\rho]} = (\rho^-)^{(\sigma g)^{-1}} (k_1[g,\omega \rho])$.
\item
(Simplification)
The identities are sorted by length and then lexicographically on the relators.
An identity is discarded if it is the empty list; if it is equal to
an earlier identity; or if it is the inverse of an earlier identity.
Furthermore, a proper subsequence is deleted if it is a conjugate of 
another identity, and the list is resorted.
\item
(Output)
The resulting list  $\cY_{\cP}$  of Y-sequences, representing 
a complete set of generators for  $\Pi_2(\cP)$  as a  $\bbZ G$-module, 
is output.
\end{enumerate}
\end{alg}

%%%%%%%%%%%%%%%%%%%%%%%%%%%%%%%%%%%%%%%%%%%%%%%%%%%%%%%%%%%%%%%%%%%%%%%%%%%%

\section{Examples} 

In this section we consider three examples.
First we return to the presentation  $\cQ$  of  $Q_8$  and the 
complete logged rewrite system obtained in Example \ref{exquat}.
We show the results obtained using the {\GAP} implementation,
obtaining  $32$  identities, which can be reduced to $18$.
Secondly we consider the free abelian group on two generators
and verify that all the identities are trivial.
Finally we consider the infinite  
$2$-generator, $1$-relator \emph{trefoil} group
and obtain a logged rewrite system which is complete with respect to a wreath
product order.

\begin{example}
\emph{The Cayley graph  $\wtX$  of the quaternion presentation
$\cQ$  with:
}
$$
X = \{a,b\}, \;\; 
R = \{\rho_1, \rho_2, \rho_3, \rho_4\}, \;\;
\omega : \rho_1 \mapsto a^4, \;       \rho_2 \mapsto b^4, \;
         \rho_3 \mapsto abab^{-1}, \; \rho_4 \mapsto a^2b^2 ~.
$$
\emph{is shown below, where the elements of  $Q_8$  are taken
to be the irreducible words obtained in Example \ref{exquat}:
}
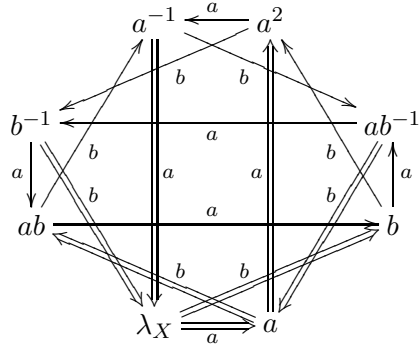
$$
\xymatrix{
   & a^{-1} \ar@{=>}[ddd]^a \ar[drr]_(0.4)b  
     & a^2 \ar[l]_a \ar[dll]^(0.4)b 
       &  \\  
 b^{-1} \ar[d]_a \ar@{=>}[ddr]^(0.4)b 
   & & &  ab^{-1} \ar[lll]^a \ar@{=>}[ddl]_(0.4)b \\  
 ab  \ar[rrr]^a \ar[uur]_(0.4)b
   & & &   b  \ar[u]_a \ar[uul]^(0.4)b \\  
   & \idX \ar@{=>}[r]_a \ar@{=>}[urr]^(0.4)b  
     & a  \ar@{=>}[uuu]^a \ar@{=>}[ull]_(0.4)b
       &  \\ }  
$$
\emph{The edges $[g,x]$ for which 
$N((\sigma g)x) = (\sigma g)x$ 
are marked with double lines in the graph. 
The image of these edges under $k_1$ is the identity. 
The images of the other edges under $k_1$, calculated using the logged
complete rewrite system $\cL(\cQ)$, 
are shown in the Table 2.  %%\ref{k1values}
}

\begin{table}[!htp] \label{k1values}
\begin{center}
\begin{tabular}{|r|c|c|l|}
\hline   
$[g,x]$  &  $w = \sigma(g(\phi x))$
  &  $\mu( (\sigma g) x w^{-1})$  &  $k_1[g,x]$\\ 
\hline 
$[a^{-1},b]$  & $ab^{-1}$ & $a^-b^+b^+a^-$    
  & $(\rho_1^-) \,(\rho_4^+)^{a^-}$  \\  
$[b,a]$       & $ab^{-1}$ & $b^+a^+b^+a^-$         
  & $(\rho_3^+)^{a^+} \,(\rho_1^-) \,(\rho_4^+)^{a^-}$ \\  
$[b,b]$       & $a^2$     & $b^+b^+a^-a^-$      
  & $(\rho_4^-) \,(\rho_2^+)^{a^{-2}}$  \\  
$[b^{-1},a]$  & $ab$      & $b^-a^+b^-a^-$ 
  & $(\rho_4^-) \,(\rho_3^+)^{a^-}$  \\
$[a^2,a]$     & $a^{-1}$  & $a^+a^+a^+a^+$
  & $(\rho_1^+)$  \\  
$[a^2,b]$     & $b^{-1}$  & $a^+a^+b^+b^+$ 
  & $(\rho_4^+)$  \\  
$[ab,a]$      & $b$       & $a^+b^+a^+b^-$   
  & $(\rho_3^+)$  \\  
$[ab,b]$      & $a^{-1}$  & $a^+b^+b^+a^+$
  & $(\rho_4^-)^{a^-} \,(\rho_2^+)^{a^{-3}} \,(\rho_1^+)$ \\  
$[ab^{-1},a]$ & $b^{-1}$  & $a^+b^-a^+b^+$
  & $(\rho_4^-)^{a^-} \,(\rho_3^+)^{a^{-2}} \,(\rho_4^+)$ \\
\hline
\end{tabular}
\end{center}
\caption{Chosen values for  $k_1[g,x]$}
\end{table}

\begin{table}[!hbp] \label{identslist}
\begin{center}
\begin{tabular}{|l|l|}
\hline
\emph{cycle} & \emph{identity} \\
\hline 
$[\idX, \rho_3]$ &
$\iota_{1} =
(\rho_1^-) \,(\rho_1^+)^{a^+}$\\
$[a^2,\rho_1]$ &
$\iota_{2} =
(\rho_1^-) \,(\rho_1^+)^{a^{+2}}$\\
$[\idX,\rho_2]$ &
$\iota_{3} =
(\rho_2^-) \,(\rho_4^-) \,(\rho_2^+)^{a^{-2}} \,(\rho_4^+)$\\
$[b^{-1},\rho_2]$ &
$\iota_{4} =
(\rho_2^-) \,(\rho_4^-)^{b^-} \,(\rho_2^+)^{a^{-2}b^-} \,
 (\rho_4^+)^{b^-}$\\
$[b,\rho_2]$ &
$\iota_{5} =
(\rho_2^-) \,(\rho_4^-)^{b^+} \,(\rho_2^+)^{a^{-2}b^+} \,(\rho_4^+)^{b^+}$\\
$[ab^{-1},\rho_3]$ &
$\iota_{6} =
(\rho_3^-) \,(\rho_4^-)^{b^-} \,(\rho_3^+)^{a^-b^-} \,
 (\rho_4^+)^{a^+b^-}$\\
$[a^2,\rho_4]$ &
$\iota_{7} =
(\rho_4^-) \,(\rho_1^+)^{a^{+2}} \,(\rho_4^-)^{a^{+2}} \,(\rho_2^+)$\\
$[a^{-1},\rho_4]$ &
$\iota_{8} =
(\rho_4^-) \,(\rho_4^-)^{a^{-2}} \,(\rho_2^+)^{a^{-4}} \,
 (\rho_1^+)^{a^-}$\\
$[b,\rho_b]$ &
$\iota_{9} =
(\rho_4^-) \,(\rho_3^+)^{a^+b^+} \,(\rho_1^-)^{b^+} \,
 (\rho_3^+)^{a^{-2}b^+} \,(\rho_4^+)^{b^+}$\\
$[ab,\rho_4]$ &
$\iota_{10} =
(\rho_4^-) \,(\rho_3^+)^{a^+b^+} \,(\rho_3^+)^{a^{+2}b^+} \,
 (\rho_1^-)^{a^+b^+} \,(\rho_4^+)^{b^+}$\\
$[b,\rho_1]$ &
$\iota_{11} =
(\rho_1^-) \,(\rho_3^+)^{a^+b^+} \,(\rho_1^-)^{b^+} \,(\rho_3^+)^{a^{-2}b^+} \,
 (\rho_3^+)^{a^-b^+} \,(\rho_3^+)^{b^+}$\\
$[ab,\rho_1]$ &
$\iota_{12} =
(\rho_1^-) \,(\rho_3^+)^{a^+b^+} \,(\rho_3^+)^{a^{+2}b^+} \,(\rho_1^-)^{a^+b^+}
 \,(\rho_3^+)^{a^-b^+} \,(\rho_3^+)^{b^+}$\\
$[b,\rho_3]$ &
$\iota_{13} =
(\rho_3^-) \,(\rho_3^+)^{a^+b^+} \,(\rho_1^-)^{b^+} \,(\rho_4^+)^{a^-b^+} \,
 (\rho_2^-)^{a^{-2}b^+} \,(\rho_4^+)^{b^+}$\\
$[ab,\rho_3]$ &
$\iota_{14} =
(\rho_3^-) \,(\rho_3^+)^{a^+b^+} \,(\rho_4^-)^{a^+b^+} \,
 (\rho_2^+)^{a^-b^+} \,(\rho_2^-)^{a^{-2}b^+} \,(\rho_4^+)^{b^+}$\\
$[b^{-1},\rho_3]$ &
$\iota_{15} =
(\rho_3^-) \,(\rho_4^-)^{b^-} \,(\rho_3^+)^{a^-b^-} \,
 (\rho_4^-)^{a^-b^-} \,(\rho_2^+)^{a^{-3}b^-} \,(\rho_1^+)^{b^-}$\\
$[ab^{-1},\rho_4]$ &
$\iota_{16} =
(\rho_4^-) \,(\rho_4^-)^{b^-} \,(\rho_3^+)^{a^-b^-} \,(\rho_3^+)^{b^-} \,
 (\rho_4^-)^{b^-} \,(\rho_2^+)^{a^{-2}b^-} \,(\rho_4^+)^{b^-}$\\
$[b^{-1},\rho_1]$ &
$\iota_{17} =
(\rho_1^-) \,(\rho_4^-)^{b^-} \,(\rho_3^+)^{a^-b^-} \,
 (\rho_3^+)^{b^-} \,(\rho_3^+)^{a^+b^-} \,(\rho_1^-)^{b^-} \,
 (\rho_3^+)^{a^{-2}b^-} \,(\rho_4^+)^{b^-}$\\
$[ab^{-1},\rho_1]$ &
$\iota_{18} =
(\rho_1^-) \,(\rho_4^-)^{b^-} \,(\rho_3^+)^{a^-b^-} \,
 (\rho_3^+)^{b^-} \,(\rho_3^+)^{a^+b^-} \,(\rho_3^+)^{a^{+2}b^-} \,
 (\rho_1^-)^{a^+b^-} \,(\rho_4^+)^{b^-}$\\
\hline
\end{tabular}
\end{center}
\caption{Non-trivial identities  $\iota_{[g,\rho]}$}
\end{table}

\emph{Each of the  $32$  relator cycles  $[g,\rho]$ 
is split into its component edges in  $\wtX$.
Where the direction of the cycle is contrary to the direction
of the generator on an edge, the inverse edge is used. 
We thus obtain the identity 
}
$$
\iota_{[g,\rho]} \;=\; (\rho^-)\,(k_1[g,\rho])^g \,.
$$
\emph{For example: 
}
$$
[ab,\rho_3] \mapsto [ab,a][b,b][a^2,a][a^{-1},b^{-1}]
  \, = \, [ab,a][b,b][a^2,a][ab,b]^{-1} ~,
$$
\emph{and $k_1$ of this product is read off from the chosen values for $k_1$ 
on the four edges. 
The resulting list of  $32$  identities may be shortened by omitting
trivial identities and duplicates.  
The identity
}
$$
\iota_{[a^{-1},\,\rho_2]} =
(\rho_2^-) \,(\rho_1^-)^{a^-} \, (\rho_2^+)^{a^{-4}} \,(\rho_1^+)^{a^-}
$$
\emph{satisfies the primary identity property (see Remark \ref{primidy})
and so may be omitted.
The remaining $18$ identities, ordered by length,
are given in Table 3.  %% \ref{identslist}.
The list represents a complete set of generators for the  
$\bbZ G$-module of identities among relations 
for the presentation  $\cQ$  of  $Q_8$.
In fact this set can be reduced to  $6$  generators, 
but the reduction requires methods dealt with in \cite{HeRe}.
}
\end{example}

\begin{example}
\emph{Our second example is the infinite abelian group with presentation
$\cA = \grp \langle X, \, \omega : R \to F(X) \rangle$ 
where  $X = \{x,y\}$
and  $R = \{\rho\}$  with   $\omega(\rho) = xyx^{-1}y^{-1}$.
The initial logged rewrite system is:
\begin{align*}
\cL_{init} \,=\, & \{~(x^+x^-, \idbY, \idbX), \;
                      (x^-x^+, \idbY, \idbX), \;
                      (y^+y^-, \idbY, \idbX), \;
                      (y^-y^+, \idbY, \idbX), \\
               & \;\; (x^+y^+x^-y^-, (\rho^+), \idbX )~\}. 
\end{align*}
Logged Knuth-Bendix completion terminates, yielding the system
\begin{align*}
\cL(\cA) \,=\, &  \{~
      (x^+x^-, \idbY, \idbX), \; (x^-x^+, \idbY, \idbX), \;
      (y^+y^-, \idbY, \idbX), \; (y^-y^+, \idbY, \idbX), \\
           & \;\;
      (y^+x^+,(\rho^-),x^+y^+), \;
      (y^+x^-,(\rho^+)^x,x^-y^+), \\
           & \;\;
      (y^-x^+,(\rho^+)^{xyx^{-1}},x^+y^-), \;
      (y^-x^-,(\rho^-)^{xy},x^-y^-)~\}. 
\end{align*}
It can be deduced that the set of elements of $\bar{X}$ 
which are irreducible with respect to the complete rewrite system  $\cL(\cA)$  
is  $\{x^n y^m : n,m \in \bbZ \}$.
This enables us to say that the (infinite) Cayley graph $\wtX$ has
edges of the form  $[x^n y^m, x]$  and  $[x^n y^m, y]$. 
We must now define  $k_1$  on all such edges.
First note that there are two cases when  $k_1$  maps an edge to  $\idbY$, 
firstly when the edge is of the form  $[x^n y^m, y]$, 
since  $x^n y^{m+1}$  is irreducible, 
and secondly when the edge is of the form  $[x^n, x]$.
It remains to determine  $k_1[x^n y^m, x]$  when  $m>0$  and when  $m<0$. 
Using the formula for  $k_1$  in the proof of Theorem \ref{sepdef2},
}
$$
k_1[x^n y^m, x]
 = L(x^n y^m x \, N(x^n y^m x)^{-1}) )
 = L(x^n y^m x y^{-m} x^{-(n+1)}) ~.
$$
\emph{When  $m>0$,  logged rewriting of 
$x^n y^m x y^{-m} x^{-(n+1)}$  gives 
}
$$ k_1[x^n y^m, x] =
  (\rho^-)^{y^{-(m-1)}x^{-n}} \,(\rho^-)^{y^{-(m-2)}x^{-n}} \cdots
      (\rho^-)^{y^{-1}x^{-n}} \,(\rho^-)^{x^{-n}}
$$ 
\emph{by repeated application of the logged rule  
$yx \rightarrow xy \; \mbox{by} \; (\rho^-)$.
For  $m<0$,  repeated application of 
$y^{-1}x \rightarrow xy^{-1} \; \mbox{by} \; (\rho^+)^{xyx^{-1}}$  gives
\begin{eqnarray*}
  &   &  k_1[x^n y^m, x] \\
  & = & 
  (\rho^+)^{xyx^{-1}y^{-m-1}x^{-n}} \,(\rho^+)^{xyx^{-1}y^{-m}x^{-n}} \cdots
    (\rho^+)^{xyx^{-1}y^{-2}x^{-n}} \,(\rho^+)^{xyx^{-1}y^{-1}x^{-n}} ~.
\end{eqnarray*}
There are infinitely many pairs  $[g,\rho]$  for this group, 
but each has the form  $[x^n y^m, \rho]$. 
The boundary of the relator cycle is  
}
$$
[x^n y^m, xyx^{-1}y^{-1}] \; = \;
[x^ny^m,x]~[x^{n+1}y^m,y]~[x^n,y^{m+1},x]^{-1}~[x^ny^m,y]^{-1} ~.
$$
\emph{The image under  $k_1$  of this cycle is 
$(\rho^+)^{y^{-m}x^{-n}}$  when  $m>0$,  and
$(\rho^+)^{xyx^{-1}y^{-m-1}x^{-n}}$  when  $m<0$. 
The identities formula in Theorem \ref{sepdef1} for 
$\sep[g,\rho]$  thus gives us 
}
$$
(\rho^-)^{y^{-m}x^{-n}} \,(\rho^+)^{(x^ny^m)^{-1}} \,=_P\, \idbY, 
\quad
(\rho^-)^{xyx^{-1}y^{-m-1}x^{-n}} \,(\rho^+)^{(x^ny^m)^{-1}} \,=_P\, \idbY ~.
$$
\emph{So we have verified algebraically that all the identities among
relations for this one relator group are trivial. 
}
\end{example}

\begin{example} \label{trefoil}
\emph{The application of logged rewrite systems to the work of 
\cite{BrSa} allows the direct computation of a finite set of generators 
for the module of identities among relations for presentation 
of a finite group. 
In the case of infinite groups the
computation is more difficult as it relies upon the 
logged reduction of general expressions such as 
$x^n y^m x y^{-m} x^{-(n+1)}$.
An infinite group where the problem of proving asphericity algebraically
using the Brown-Razak set of generators for the module of identities
among relations 
is not straightforward is the trefoil knot group, with presentation
}
$$
\cT = \grp \langle X = \{x,y\}, \;
\omega : R = \{\rho\} \to F(X), \, \rho \mapsto x^3 y^{-2} \rangle ~.
$$
\emph{In order to obtain a complete rewrite system for  $\cT$
it is convenient to use a wreath product or \emph{syllable} ordering 
(see \cite{LeCh}).
If  $X = \{x_1, x_2, \ldots, x_n\}$  first choose
$x_1^- > x_1^+ > x_2^- > \cdots > x_n^+$.
For  $S \subset \bar{X}$  denote by  $\max S$  the largest generator in  $S$.
For each  $u \in S^*$  denote by  $n_S(u)$  the number of occurrences
of  $\max S$  in  $u$.
The order on  $\bar{X}^*$  is defined recursively, taking
}
$$
S \; \in \; \{ ~ \{x_n^+\}, ~ \{x_n^-, x_n^+\}, ~\ldots,
               ~ \bar{X} \setminus \{x_1^-\}, ~ \bar{X} ~\}.
$$
\emph{For  $u,v \in S^*$  we take  $u < v$  provided:
}\begin{itemize}
\item $n_S(u) < n_S(v)$, \emph{or}
\item $n_S(u) = n_S(v) = n, \; m = \max S$, \emph{and there exists}
      $0 \le i \le n$  \emph{and}  \\
      $w_1, \ldots, w_i, u_{i+1}, \ldots, u_{n+1}, v_{i+1}, \ldots, v_{n+1}
        \in (S \setminus \{m\})^*$  \emph{with}  $u_{i+1} < v_{i+1}$.
\end{itemize}
\emph{The complete rewrite system for  $\cT$  with respect to this 
syllable ordering has been obtained using the COSY completion system
at Kaiserslautern (see \cite{MaRe}), and the logged information
was later recovered by hand, resulting in the following 
complete logged system:
}\begin{eqnarray*}
\cL(\cT)
 & = &  \{ (y^+y^-, \idbY, \idbX),\,
           (y^-y^+, \idbY, \idbX),\,
           (x^{+3}, (\rho^+), y^{+2}), \\
 &     &   (y^{+2}x^+, (\rho^-)(\rho^+)^{x^{-1}}, x^+y^{+2}),\, 
           (y^-x^+, (\rho^-)^{x^{-1}y}(\rho^+)^y, y^+x^+y^{-2}), \\
 &     &   (x^-, (\rho^-)^x, x^{+2}y^{-2}) \}.
\end{eqnarray*}

\emph{Although the normal forms of the group elements are too 
irregular to attempt to describe the Brown-Razak identities, 
we can consider the identities that arise from the logged reduction 
of the critical pairs of  $\cL(\cT)$. 
There are six words on which the rules overlap, and we expect that 
the six identities which arise from these are trivial.
}

\emph{The six overlap words are:
}
$$
y^+y^-y^+, \;  y^-y^+y^-, \;  y^-y^{+2}x^+, \; 
y^+y^-x^+, \;  y^-x^{+3} \; \mbox{ and } \; y^{+2}x^{+3}.
$$
\emph{The logged reductions of the first two immediately give the trivial 
$Y$-sequence as the identity, so we consider the word  $y^-y^{+2}x^+$. 
It can be reduced to 
$(\idbY,y^+x^+)$ or to}\\
$((\rho^-)^y(\rho^+)^{x^{-1}y}(\rho^-)^{x^{-1}y}(\rho^+)^y,y^+x^+)$.
\emph{Thus the identity here is  
}
$$
(\rho^-)^y(\rho^+)^{x^{-1}y}(\rho^-)^{x^{-1}y}(\rho^+)^y,
$$ 
\emph{which is equivalent under  $=_P$  (Peiffer and inverse rules) 
to the identity $Y$-sequence $\idbY$. 
Similarly,  $y^+y^-x^+$  can be reduced to  $(\idbY,x^+)$ 
and 
$((\rho^-)^{x^{-1}} (\rho^+) (\rho^-) (\rho^+)^{x^{-1}},x^+)$, 
giving the identity
}
$$
(\rho^-)^{x^{-1}} (\rho^+) (\rho^-) (\rho^+)^{x^{-1}}=_P \idbY.
$$
\emph{The fifth overlap word is  $y^-x^{+3}$, 
and can be reduced to  $((\rho^+)^y,y^+)$  and to 
}
$$
((\rho^-)^{x^{-1}y} (\rho^+)^y (\rho^-)^{x^{-1}y^2x^{-1}y^{-1}}
(\rho^+)^{y^2x^{-1}y^{-1}} (\rho^+)^{y^{-1}} (\rho^-)^{x^{-1}y^2xy^{-3}}
(\rho^+)^{y^2xy^{-3}},y^+)~, 
$$
\emph{giving the identity
}
$$
(\rho^-)^{x^{-1}y} (\rho^+)^y (\rho^-)^{x^{-1}y^2x^{-1}y^{-1}}
(\rho^+)^{y^2x^{-1}y^{-1}} (\rho^+)^{y^{-1}} (\rho^-)^{x^{-1}y^2xy^{-3}}
(\rho^+)^{y^2xy^{-3}}(\rho^-)^y =_P \idbY.
$$ 
\emph{The final overlap word  $y^{+2}x^{+3}$  gives the identity
}
$$
(\rho^-) (\rho^+)^{x^{-1}} (\rho^-)^{x^{-1}} (\rho^+)^{x^{-2}} 
(\rho^-)^{x^{-2}} (\rho^+) (\rho^+)^{y^{-2}} (\rho^-)^{y^{-2}} =_P \idbY.
$$

\emph{This example shows how to obtain identities among relations from critical
pairs of a complete logged rewrite system. 
The paper \cite{BHW} will pursue this further 
and relate it to work on derivation schemes.
}
\end{example}

%%%%%%%%%%%%%%%%%%%%%%%%%%%%%%%%%%%%%%%%%%%%%%%%%%%%%%%%%%%%%%%%%%%%%%%%%%%

\end{document}